\documentstyle[12pt]{article}

\parskip=\medskipamount                 
\def\beee{\begin{equation}}
\def\eeee{\end{equation}}

\def\beq{\begin{eqnarray}}

\def\eeq{\end{eqnarray}}

\begin{document}

\rightline{Preprint CRN 94/08 and IFUNAM FT 94-41}
\rightline{February 1994}
\rightline{May 1994 (corrected)}
\rightline{funct-an/9403002}

\begin{center}
{\Large \bf On Polynomial  Relations in the Heisenberg Algebra}\\

\vskip 2cm

Norbert Fleury \\
CRN , Physique Th\'eorique and  Universit\'e Louis Pasteur,
B.P. 20 F-67037 Strasbourg Cedex \\

and\\

Alexander Turbiner
\footnote{On leave in absence from Institute for Theoretical and Experimental
Physics,
Moscow, Russia\\ e-mail: turbiner@teorica0.ifisicacu.unam.mx}
\\
Instituto de Fisica, UNAM, Apartado Postal 20-364, 01000 Mexico D.F., Mexico\\

\vskip 2.5cm

{\bf Abstract}\\
\begin{quote}
Polynomial relations between the generators of the classical and quantum
Heisenberg
algebras are presented. Some of those relations can have a meaning of the
formulas of the
normal ordering  for the creation/annihilation operators occurred in the method
of the second quantization.
\end{quote}
\end{center}

\newpage
\indent
Undoubtedly, the Heisenberg algebra plays an exceptionally important role in
many
branches of the theoretical and mathematical physics. In particular,  the
Heisenberg
algebra appears as a basic element of the method of the second quantization
(see, e.g. \cite{1} ) in quantum mechanics, quantum field theory, statistical
mechanics
and also in nuclear physics, in the interacting boson model
by Arima-Iachello \cite{2}.  Special importance of the Heisenberg algebra stems
from
the fact, that any simple or semi-simple Lie algebra can be found
as an embedding in the universal enveloping algebra of the Heisenberg algebra.

This paper is mostly  devoted to the investigation of non-linear (polynomial)
relations
between the elements of the universal enveloping algebra of the Heisenberg
algebra  as a consequence either existence of finite-dimensional
representations of the classical Lie
algebras or presence of two-dimensional Borel subalgebra in
any simple or semi-simple Lie algebra.
Surprisingly, all relations we found can be extended to the case of recently
discovered
quantum ($q-$deformed) algebras. Some of those relations have a meaning
of the formulas of the normal (lexicographical) ordering. An existence of those
relations
simplifies drastically a solution one of the the most tedious (and most common)
problem appearing in concrete calculations in the framework of the second
quantization method
-- the problem of the ordering of polynomials in creation/annihilation
operators.
It is worth noting that  the problem of ordering also appears in representation
theory
of finite- and infinite-dimensional Lie algebras. Above-mentioned relations
play
an important role in a problem of a classification of recently-discovered {\it
quasi-exactly-solvable} Schroedinger equations -- the eigenvalue problem
for the Schroedinger operators, where a certain
amount of eigenstates can be found algebraically (see e.g. \cite{8}). In
particular, they
simplify finding a number of free parameters of the quasi-exactly-solvable
Schroedinger
operators.

\noindent
1. Take two linear operators $a, b$  obeying the
\begin{equation}
\label{e1}
[ a, b ] =1
\end{equation}
Under a certain conditions those operators can be interpreted as
creation/annihilation operators.
One can easily verify by the mathematical induction that
\begin{equation}
\label{e2}
ba(ba-1)(ba-2)\ldots(ba-n)=b^{n+1}a^{n+1}
\qquad n=0,1,2,\ldots
\end{equation}
holds. One can show  that the opposite statement is also valid: if (2) holds,
then the operators $a,b$ must obey (1). Since the  multipliers in the l.h.s. of
(2)
commute, they  can be placed  in arbitrary order.
In general, changing a $c-$number in a bracket in the l.h.s. of (2) : $(ba-m)
\mapsto (ba-m-\epsilon) , \epsilon \neq 0$,
leads immediately to the appearance of the terms $b^k a^k , k=0,1, \ldots n$ in
the r.h.s.
Therefore, setting $\epsilon$ equals to zero implies fulfillment of $n$
conditions: the leading term only survives in the r.h.s. of (2), while all
other $n$ terms vanish.

Now let us consider  two concrete realizations of $a, b$.

{\it (a)} \ Assume that
\begin{equation}
\label{e3}
b=x,\;a=\partial_x\qquad \left(\partial_x\equiv{d\over dx}\right)
\end{equation}
then the relation (2) becomes \cite{3,4}
\begin{equation}
\label{e4}
x\partial_x(x\partial_x-1) \ldots(x\partial_x-n)=x^{n+1}\partial_x^{n+1}
\qquad n=0,1,2,\ldots
\end{equation}

{\it (b)} \ If
\[
b={1\over\sqrt2}(\partial_x+x), \quad a={1\over\sqrt2}(\partial_x-x)
\]
then
\begin{equation}
\label{e5}
(\partial_x^2-x^2-1)(\partial_x^2-x^2-3)\ldots
(\partial_x^2-x^2-2n-1)=(\partial_x+x)^{n+1}(\partial_x-x)^{n-1}
\qquad n=0,1,2,\ldots
\end{equation}

\noindent
2. Another class of relations occurs as a consequence of the existence of
two-dimensional
Borel sub-algebra
\begin{equation}
\label{e6}
[ A,B ]=A
\end{equation}
in $s\ell_2$ ($s\ell_{k+1}$) . Let (6) holds, then \cite{5}
\begin{equation}
\label{e7}
(BA)^{n+1}=(B)_{n+1} A^{n+1},\;(B)_{n+1}=B(B+1)\ldots (B+n) \quad ,
\qquad n=0,1,2,\ldots
\end{equation}
and also
\begin{equation}\label{e8}
(AB)^{n+1}=A^{n+1} (B)_{-n-1},\;(B)_{-n-1}=B(B-1)\ldots (B-n) \quad ,
\qquad n=0,1,2,\ldots
\end{equation}
are valid. The formulas (7)-(8) can be interpreted as the formulas of the
normal ordering:
on the r.h.s. of (7)-(8) one has the $A$-type operators are placed on the left,
$B$-type
operators on the right.

The algebra  $s\ell_2$ has a natural embedding into the universal enveloping
algebra of the Heisenberg algebra (1): the operators
\[
J^+\ =\ b^2a-2\alpha b
\]
\begin{equation}
\label{e9}
J^0\ =\ ba-\alpha
\end{equation}
\[
J^-\ =\ a
\]
obey $s\ell_2$-commutation relations
\[
[ J^\pm, J^0 ]=\mp J^\pm;\qquad [ J^+,J^-]=-2J^0
\]
for any $\alpha$
\footnote{If the operators $a,b$ are realized as in (3), then (9) becomes
well-known representation of $sl_2$ algebra in the first-order differential
operators,
derived by Sophus Lie}.
Taking in (6) $A=J^-,\;B=J^0$ at $\alpha=n$ and substituting it into  (7), one
gets
\begin{equation}
\label{e10}
(ba^2-na)^{n+1}=b^{n+1} a^{2n+2},\qquad n=0,1,2,\ldots
\end{equation}
The relation (10) can be compared with
\begin{equation}
\label{e11}
(b^2a-nb)^{n+1}=b^{2n+2} a^{n+1},
\qquad n=0,1,2,\ldots
\end{equation}
obtained in \cite{6} as a consequence of existence of finite-dimensional
representations of $s\ell_2$-algebra of the first order differential
operators. Under the assumption that the operators $a, b$ are
hermitian-conjugated:
\begin{equation}
\label{e12}
(a)^+=b
\end{equation}
the equalities (10) and (11) are also related by the hermitian conjugation.
Analogously,
using (8), (9) at $\alpha=0$ , it emerges that
\footnote{In fact the equality (13) is a particular case of a general Theorem:
\\
{\it If two operators $a,b$ obey (1), then the equality
\[
(abab \ldots a)^n = a^nb^n a^n b^n \ldots a^n \ , n=0,1,2,...
\]
holds}.\\
The proof is quite straightforward based mainly on an application of well-known
basic
equalities
\[
 [a^n,b]=na^{n-1}\ ,\  [b^n,a]=-nb^{n-1}.
\]}
\begin{equation}
\label{e13}
(aba)^{n+1}=a^{n+1} b^{n+1}a^{n+1}
\end{equation}
A remarkable property of (13) is that the equality is independent on the value
of  the constant
in the r.h.s. (1). In particular, the formula (13) is trivially valid for
commuting operators $a,b$.

\noindent
3. As one of possible generalizations, let us take $(2p+1)$-dimensional
Heisenberg
algebra
\begin{equation}
\label{e14}
[ a_i,a_j] =[ b_i,b_j] = 0 \ ,\quad
[ a_i,b_j] =\delta_{ij}\ ,\qquad i,j=1,2,\ldots p
\end{equation}
Then one can derive an immediate extension of (2) :
\begin{equation}
\label{e15}
\prod_{\ell=0}^n\left(\sum_{i=1}^p b_ia_i-\ell\right)=
\sum_{j_1+j_2+\cdots +j_p=n+1}
 C_{j_1\ldots j_p}^{n+1} b_1^{j_1}a_1^{j_1}
 b_2^{j_2}a_2^{j_2}\ldots b_p^{j_p}a_p^{j_p}
\end{equation}
and also a generalization of (11) \cite{6}
\[
\left\lbrack b_\ell\left(\sum_{m=1}^p b_m a_m-n\right)
\right\rbrack^{n+1} \ =
\]
\begin{equation}\label{e16}
=\ (b_\ell)^{n+1}\sum_{j_1+j_2+\cdots j_p=n+1}
 C_{j_1j_2\ldots j_p}^{n+1}
(b_1)^{j_1}(b_2)^{j_2}\ldots(b_p)^{j_p}
a_1^{j_1}a_2^{j_2}\ldots a_p^{j_p}
\end{equation}
and (10)
\[
\left\lbrack \left(\sum_{m=1}^p b_ma_m-n\right)a_\ell
\right\rbrack^{n+1}\ =
\]
\begin{equation}
\label{e17}
=\ \sum_{j_1+j_2+\cdots j_p=n+1}
 C_{j_1j_2\ldots j_p}^{n+1}
(b_1)^{j_1}(b_2)^{j_2}\ldots(b_p)^{j_p}
a_1^{j_1}a_2^{j_2}\ldots a_p^{j_p} \ (a_\ell)^{n+1}
\end{equation}
where $\ell=1,2,\ldots p;\; n=0,1,2\ldots$ and  $C_{j_1j_2\ldots j_p}^{n+1}$
are the multinomial coefficients.

\noindent
4. Now let us consider $q$-deformed Heisenberg algebra
\begin{equation}
\label{e18}
[ a,b]_q\equiv ab-qba=1, \qquad q\in \Re
\end{equation}
coinciding (1) at $q=1$.

Using the mathematical induction, one can prove that
\begin{equation}
\label{e19}
ba(ba-1)(ba-\{2\})\ldots (ba-\{ n\})\ =\ q^{n(n+1)\over2}b^{n+1}a^{n+1}
\end{equation}
(cf. (2)), where $\{ n\}={1-q^n\over 1-q}$ is a so-called $q$-number.
One can show  that the opposite statement is also valid: if (19) holds,
then the operators $a,b$ must obey (18).

A natural representation of (18) is
\begin{equation}
\label{e20}
a =D, \qquad  b=x
\end{equation}
where $D$ is the Jackson symbol defined as :
\[
Df(x)={f(x)-f(qx)\over x(1-q)}.
\]
Substitution of (20) to (19) leads to \cite{6}
\begin{equation}
\label{e21}
xD(xD-1)\ldots (xD-\{ n\})=
q^{n(n+1)\over2} x^{n+1}D^{n+1}
\end{equation}
There exist  a certain $q$-generalizations of (7) and (8) as well. Let
\[
AB-qBA=A
\]
 (cf. (6)), then
\begin{equation}
\label{e22}
(BA)^{n+1}=B_{n+1,q}A^{n+1}, B_{n+1,q}=B(qB+1)\ldots (q^nB+\{ n\})
\end{equation}
and
\begin{equation}
\label{e23}
(AB)^{n+1}=A^{n+1}B_{-n-1,q}, B_{-n-1,q}=B\left({1\over q}B-{1\over q}
\right)\ldots
\left({1\over q^n}B-{\{ n\}\over q^n}
\right)
\end{equation}
that can be verified by the mathematical induction.
\medskip
Also, by straightforward calculations one shows that, once (18) holds, then
the operators
\[
\hat J^+=(1-2\alpha(1-q))^{-1/2}(b^2a-2\alpha b)
\]
\begin{equation}
\label{e24}
\hat J^0=(1-2\alpha(1-q))^{-1}\,\left\lbrack\left(1+{2\alpha(q^2-q)\over1+q}
\right)ba-{2\alpha\over1+q}
\right\rbrack
\end{equation}
\[
\hat J^-=(1-2\alpha(1-q))^{-1/2} a
\]
obey $s\ell_{2q}$-commutation relations
\[
q\hat J^0\hat J^--\hat J^-\hat J^0=-\hat J^-
\]
\begin{equation}\label{e25}
\hat J^0\hat J^+-q\hat J^+\hat J^0=\hat J^+
\end{equation}
\[
q^2\hat J^+\hat J^--\hat J^-\hat J^+=-(q+1)\hat J^0
\]
for any value of $\alpha$ at fixed $q$.
In the particular representation (20) , the generators (24) coincide with those
found in \cite{7}
(see also \cite{8}).

Take two operators
\[
A=a
\]
\[
B={ba-2\alpha\over1-2\alpha(1-q)}
\]
(cf.(24)), which obey  $q$-deformed commutation relation $AB-qBA=A$.
Choosing
$2\alpha=\lbrace n\rbrace$ and using (22) together with (19), one arrives at
\begin{equation}
\label{e26}
(ba^2-\{ n\} a)^{n+1}=q^{n(n+1)} b^{n+1} a^{2n+2}
\end{equation}
(cf.~(10)). In the particular representation (20) of $a,b$, it becomes
\begin{equation}
\label{e27}
(xD^2-\{ n\} D)^{n+1}=q^{n(n+1)}x^{n+1}D^{2n+2}
\end{equation}
It is worth noting that once (18) holds,
the following relation  \cite{6} also
\begin{equation}
\label{e28}
(b^2a-\{ n\} b)^{n+1}=q^{n(n+1)}b^{2n+2}a^{n+1}
\end{equation}
(cf. (11)) holds.

It can be shown that surprisingly the equality (13) remains valid under the
deformation
(18) of the Heisenberg algebra, if the parameter of deformation \linebreak $q
\neq 0$
\footnote{
It is worth noting that $q$-analogs of the basic equalities (see footnote 3)
are
\[
a^n b - q^n b a^n = \{ n\} a^{n-1}\ , \ b^n a - {1 \over q^n} a b^n = - {\{ n\}
\over q^n} b^{n-1}
\]
}.

\noindent
5. An attempt to generalize the results of the Section $4$ to the case of two
pairs of  the operators  $b_{1,2},\;a_{1,2}$ ( as has been done in the Section
3) demands
the implementation a certain deformation of the 5-dimensional Heisenberg
algebra
\[
 b_1b_2=\sqrt q b_2b_1,\qquad\sqrt qa_1a_2=a_2a_1
\]
\[
a_1b_1-qb_1a_1=1+(q-1)b_2a_2
\]
\[
a_2b_2-qb_2a_2=1,\qquad a_1b_2-\sqrt qb_2a_1=0
\]
\begin{equation}
\label{e29}
a_2b_1-\sqrt qb_1a_2=0
\end{equation}
in order to get needed result. Using the  mathematical induction, one can prove
that
\[
\prod_{\ell=0}^n\left(b_1a_1+b_2a_2-\{\ell \}\right)\ =
\]
\begin{equation}
\label{e30}
=\ \sum_{\ell=0}^{n+1} q^{{1\over2}n(n+1)+{\ell\over2}(\ell-n-1)}
{n+1\choose\ell}_q (b_1)^{n+1-\ell}(b_2)^\ell
 a_1^{n+1-\ell}a_2^\ell\ ,\ n=0,1,2,
\end{equation}
(cf.~(15) at $p=2$) where
${n \choose k}_q\equiv {\{ n\}! \over \{ k\}!\{ n-k\}!}\ ,\ \{ n\}= \{ 1\}
\{2\}\cdots\{ n\}$
are $q$-binomial coefficient and $q$-factorial, respectively.
\medskip
One can easily show that $q$-commutation relation
\begin{equation}
\label{e31}
\lbrack a_1\ ,\ (b_1a_1+b_2a_2-\alpha)
\rbrack_q\ =\ (1-\alpha+\alpha q) a_1\ ,
\end{equation}
is valid for any value of the parameter $\alpha$. Then denote: $A=a_1$,
$B={b_1a_1+b_2a_2-\alpha \over 1-\alpha +\alpha q}$, take $\alpha=\{n\}$ and
plug them
into (22). This leads to
\[
(b_1a_1^2+ b_2a_2a_1-\{ n\} a_1)^{n+1}\ =
\]
\begin{equation}
\label{e32}
=\ \sum_{\ell=0}^{n+1} q^{n(n+1)+{\ell \over 2}(\ell-n-1)}
{n+1\choose\ell}_q (b_1)^{n+1-\ell}(b_2)^\ell
 a_1^{n+1-\ell}a_2^\ell a_1^{n+1}
\end{equation}

If the $q$-deformed commutator $[A,B]_q=B$, then one has $(BA)^{n+1}=B^{n+1}
A_{n+1,q}$. Applying this result to
$A={b_1a_1+b_2a_2-\alpha\over 1-\alpha +\alpha q}$, $B=b_1$, and choosing
$\alpha=\{n\}$, one gets the slightly different ordering formula
\[
(b_1^2a_1+ b_1b_2a_2-\{ n\} b_1)^{n+1}\ =
\]
\begin{equation}
\label{e33}
=\ \sum_{\ell=0}^{n=1} q^{n(n+1)+{\ell \over 2}(\ell-n-1)}{n+1\choose\ell}_q
 (b_1)^{2n+2-\ell}(b_2)^\ell a_1^{n+1-\ell}a_2^\ell
\end{equation}
which turns out to be the hermitian conjugated of  (32).

Take a particular realization of the algebra (29)
\begin{equation}
\label{e34}
b_1=x,\;b_2=y,\;a_1=D_x,\;a_2=D_y
\end{equation}
then the algebraic relations (29) become the rules of $q$-calculus of the
quantum
plane introduced by Wess and Zumino \cite{9}. A substitution of (34) to (30),
(32)-(33) leads to some operator identities for finite-difference operators,
in particular, for the case of (33) they coincide to those described in
\cite{6}.
\vskip 1.5truecm
One of us (A.T.) wants to thank the Centre de Recherche Nucl\^{a}ire, Physique
Th\'eorique, Case Western Reserve University and Instituto de Fisica, UNAM,
where
this work has been done,  for kind hospitality extended to him. Also A.T. is
very grateful
to Prof. R.~Askey for valuable discussions and bringing attention to the
paper \cite{5} and Prof. S.~Szarek for an interesting suggestion.
A.T. is supported in part by CAST grant of US National Academy of Sciences.

\end{document}